\documentclass{amsart}
\usepackage{amsmath,amssymb}
\usepackage{url,hyperref}

\input xy
\xyoption{all}

\newtheorem{theorem}{Theorem}[section]
\newtheorem*{theorem*}{Theorem}
\newtheorem{lemma}[theorem]{Lemma}

\newtheorem{corollary}[theorem]{Corollary}
\newtheorem{proposition}[theorem]{Proposition}

\theoremstyle{definition}

\newtheorem{example}[theorem]{Example}
\newtheorem{remark}[theorem]{Remark}

\def\Hilb{{\operatorname{Hilb}}}
\def\Sym{{\operatorname{Sym}}} 
\def\rank{{r}}
\def\<{{\langle}}
\def\>{{\rangle}}
\def\corank{{cr}}

\title{Products of Linear Forms\\{\tiny and}\\Tutte Polynomials}
\author{Andrew Berget}
\address{Department of Mathematics,
One Shields Avenue,
University of California,
Davis, CA 95616
}
\email{berget@math.ucdavis.edu}

\date{\today}

\begin{document}
\pagestyle{plain}
\begin{abstract}
  Let $\Delta$ be a finite sequence of $n$ vectors from a vector space
  over any field. We consider the subspace of $\operatorname{Sym}(V)$
  spanned by $\prod_{v \in S} v$, where $S$ is a subsequence of
  $\Delta$. A result of Orlik and Terao provides a doubly indexed
  direct sum decomposition of this space. The main theorem is that the
  resulting Hilbert series is the Tutte polynomial evaluation
  $T(\Delta;1+x,y)$.  Results of Ardila and Postnikov, Orlik and
  Terao, Terao, and Wagner are obtained as corollaries.
\end{abstract}

\maketitle

\section{Introduction and Statement of the Theorem}
Let $V$ be a vector space of dimension $\ell$ over a field $K$ of
arbitrary characteristic. Let $\Delta = (\alpha_1,\dots,\alpha_n)$ be
a sequence of elements spanning $V^*$, the dual space of $V$. We allow
the possibility that some of $\alpha_i$ are the zero form, and for
some $\alpha_i$'s to differ by scalars. This is to say, $\Delta$ is a
realization of a rank $\ell$ matroid $M(\Delta)$ on ground set $E =
\{1,2,\dots,n\}$. We assume that the reader is familiar with the
basics of matroids and their Tutte polynomials (see
\cite{tutte-poly,matroid-th}).

Denote the symmetric algebra of $V^*$ by $\Sym(V^*)$, which is thought
of as the $K$-algebra of polynomial functions on $V$. Let $P(\Delta)$
be the $K$-subspace of $\Sym(V^*)$ spanned by $\alpha_S =\prod_{i \in
  S} \alpha_i$, where $S \subset E$. Let $P(\Delta)_{j,k}$ be the
$K$-span of those products $ \alpha_S$ where $\{\alpha_i:i \in E- S\}$
spans a $j$-dimensional subspace of $V^*$ and $k=|E-S|$.  It follows
from a result of Orlik and Terao (see \cite{orlik-terao,terao} and
{Proposition~\ref{prop:directsum}}) that there is a $K$-vector space
decomposition
\begin{align}\label{eq:directsum}
  P(\Delta)=\bigoplus_{0 \leq j \leq k \leq n } P(\Delta)_{j,k}.
\end{align}
The main result of this paper is the following.
\begin{theorem}\label{thm:main}
  The Tutte polynomial of $M(\Delta)$ is equal to
  \begin{align*}
    \sum_{0 \leq j \leq k \leq n} (x-1)^{\ell-j} y^{k-j} \dim
    P(\Delta)_{j,k}.
  \end{align*}
\end{theorem}
This result was anticipated in the work of many authors. Following a
suggestion of Aomoto, in \cite{orlik-terao} Orlik and Terao considered
a vector space related to $P(\Delta)$.  They studied an algebra whose
underlying $K$-vector space was isomorphic to $\bigoplus_{j = 0}^n
P(\Delta)_{j,j}$. In \cite{wagner}, Wagner considered an algebra whose
underlying $K$-vector space was isomorphic to $P(\Delta)$. In
\cite{ardila-postnikov} Ardila and Postnikov investigate the spaces
$P(\Delta)$ and $\bigoplus_{k = 0}^n P(\Delta)_{l,k}$ from the point
of view of \textit{power ideals}. Other spaces related to $P(\Delta)$
were studied in
\cite{ardila,brion-vergne,cordovil,forge-las-vergnas,proudfoot-speyer,schenck-tohaneanu}
and \cite{terao}. There is also a vast literature on box splines and
their relationship to $P(\Delta)$ and its subspaces. In this area,
Dahmen and Miccelli were considering related objects as early as 1983.
A collection of relevant references to this area can be found in De
Concini and Procesi \cite{dpb}.

It is worth noting that not all algebraic invariants of a matroid or
vector configuration are specializations of the Tutte polynomial, and
hence are not related to $P(\Delta)$. Two particular objects of
interest to the author are the Whitney algebra of a matroid, defined
by Crapo, Rota and Schmitt in \cite{cs} and the smallest general
linear group representation containing a fixed decomposable tensor
\cite{berget}.

This paper is organized as follows: We start by setting up some
notation and compute an example in Section~\ref{sec:ex}. In
Section~\ref{sec:cor} we list some corollaries of
Theorem~\ref{thm:main}. Section~\ref{sec:alg} gives a proof of a
formula, due to Terao \cite{terao}, for the Hilbert series of the
algebra generated by the \textit{reciprocals} of linear forms.  The
benefit of this proof is that it works over an arbitrary field $K$,
whereas Terao assumed that the characteristic of $K$ was zero. In
Section~\ref{sec:ff} the following question is answered: For what $d$
does $P(\Delta)$ contain $\Sym^d(V^*)$? In Section~\ref{sec:proof} we
give the deletion-contraction proof of Theorem~\ref{thm:main}, while
in Section~\ref{sec:proof2} we give a short prove of the theorem using
the seemingly weaker result stated in Corollary~\ref{cor:x=0}.

\section{An Example}\label{sec:ex}
Before proceeding, we state a refinement of the decomposition
\eqref{eq:directsum}, due to Orlik and Terao \cite{orlik-terao}, and
use this to give an example of Theorem~\ref{thm:main}. To do this we
recall some notation. For a subset $S \subset E$, the dimension of the
$K$-span of $\{\alpha_i: i \in S\}$ is called its rank (in
$M(\Delta)$) and denoted $\rank(S)$. A set $S \subset E$ is said to be
\textit{independent} (in $M(\Delta)$) if $\rank(S) = |S|$ and
\textit{dependent} if $\rank(S)< |S|$. A \textit{flat of $\Delta$} is
a set $X \subset E$ such that for all strict containments $X \subset
Y$, $\rank(X) < \rank(Y)$. The collection of flats of $\Delta$,
$L(\Delta)$, is a geometric lattice where the rank of $X$ is
$\rank(X)$. The \textit{closure} or \textit{span} of a subset of $E$
is the smallest flat containing it.

For $X \in L(\Delta)$, let $P(\Delta)_{X}$ be the subspace of
$P(\Delta)$ spanned by those $\alpha_S$ where $E-S$ has closure $X$.
\begin{remark}
  When $\Delta$ does not contain the zero form, it is possible to
  avoid the somewhat backwards definition of $P(\Delta)_X$ by noting
  that $P(\Delta)_X$ is isomorphic to the vector space spanned by the
  rational functions $1/\alpha_S$ where the closure of $S$ is equal to
  $X$. This topic will be discussed further in Section~\ref{sec:alg}.
\end{remark}
\begin{proposition}[Orlik-Terao {\cite[Lemma~3.2]{orlik-terao}}]\label{prop:directsum}
  There is a $K$-vector space direct sum decomposition
  \[
  P(\Delta)= \bigoplus_{X \in L(\Delta)} P(\Delta)_X.
  \]
\end{proposition}
Since $P(\Delta)$ is spanned by homogeneous elements this sum may be
refined by degree. Denote the degree $n-k$ subspace of $P(\Delta)_X$
by $P(\Delta)_{X,k}$. There is a $K$-vector space direct sum
decomposition,
\[
P(\Delta)= \bigoplus_{\substack{X \in L(\Delta)\\k \geq 0}}
P(\Delta)_{X,k}.
\]
The decomposition \eqref{eq:directsum} is obtained from the one above
by taking the ranks of the flats of $M(\Delta)$. Note that, by
definition, if $P(\Delta)_{X,k} \neq 0$ then $\rank(X) \leq k \leq
|X|$.

\begin{example}\label{ex:fano}
  Let $K=\mathbb{F}_2$ be the field with two elements and $\Delta =
  (\alpha_1,\dots,\alpha_7)$ be the seven nonzero elements of the dual
  of $V=\mathbb{F}_2^3$. The matroid $M(\Delta)$ is known as the Fano
  matroid and it is the rank three matroid whose circuits of size
  three are the three point lines of the Figure~\ref{fig:fano}.
\begin{figure}[h]
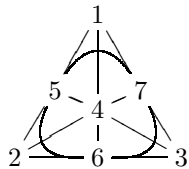

  \centering
  \[
  \xygraph{
    !{<0cm,0cm>}
    !{(0,0) }{4}="4"
    !{(0,1) }{1}="1"
    !{(0,-.5) }{6}="6"
    !{(.875,-.5) }{3}="3"
    !{(-.875,-.5) }{2}="2"
    !{(.45,.2) }{7}="7"
    !{(-.45,.2) }{5}="5"
    "1"-"5"-"2"-"6"-"3"-"7"-"1"-"4"-"6"
    "3"-"4"-"5"
    "6"-@/^15pt/"5"-@/^15pt/"7"-@/^15pt/"6"
    "2"-"4"-"7"
  }
  \]
  \caption{The Fano matroid.}
  \label{fig:fano}
\end{figure}

We will compute $T(\Delta;1+x,y)$ by finding $\dim
  P(\Delta)_{X,k}$ for all flats $X$ and all $k$ such that
  $\rank(X)\leq k \leq |X|$.
  
  See that $P(\Delta)_{\emptyset,0}$ is spanned by one nonzero element
  $\alpha_E = \prod_{i=1}^7 \alpha_i$ and hence has dimension one. The
  rank one flats of $\Delta$ are in bijection with the elements of
  $\Delta$. Hence $P(\Delta)_{\{i\},1}$ is spanned by the single
  element $\alpha_E/\alpha_i$ and $\dim P(\Delta)_{\{i\},1} = 1$.

  There are seven rank two flats of $\Delta$. If $X$ is such a flat
  then it corresponds to a set of the form
  $\{\alpha_i,\alpha_j,\alpha_i+\alpha_j\}$. It follows that
  $P(\Delta)_{X,2}$ is spanned by the three elements
  \[
  \alpha_i \alpha_{E-X},\ \alpha_j \alpha_{E-X},\
  (\alpha_i+\alpha_j)\alpha_{E-X}.
  \]
  Adding these three terms up gives $0$ and hence $\dim
  P(\Delta)_{X,2} \leq 2 $. Since none of the $\alpha_i$ are parallel,
  $\dim P(\Delta)_{X,2} = 2$. Because $P(\Delta)_{X,3}$ is spanned by
  a single nonzero element, $\dim P(\Delta)_{X,3} = 1$.

  The only rank $3$ flat of $\Delta$ is the whole set $E$. The empty
  product spans $P(\Delta)_{E,7}$ and so it has dimension one. One
  finds that $P(\Delta)_{E,6}$ is the span of
  $\alpha_1,\dots,\alpha_7$, so this space has dimension equal to the
  dimension of $V^*$, which is three. To compute $P(\Delta)_{E,5}
  \subset \Sym^2({V^*}\!)$, assume that $\alpha_1,\alpha_2$ and
  $\alpha_3$ are a basis for $V^*$. By considering leading terms
  (under any term order) we see that
  \[
  \alpha_1\alpha_2,\ \alpha_1\alpha_3,\ \alpha_2 \alpha_3,\
  \alpha_1(\alpha_1+\alpha_2),\ \alpha_2(\alpha_2+\alpha_3),\
  \alpha_3(\alpha_1 +\alpha_3)
  \]
  forms a basis for $\Sym^2({V^*}\! )$. In a similar fashion we see
  that $P(\Delta)_{E,4}$ is equal to $\Sym^3({V^*}\!)$, which has
  dimension $10$. We resort to a computer to find the dimension of
  $P(\Delta)_{E,3}$, which is spanned by $28$ products. This space is
  contained in $\Sym^4({V^*}\!)$ which has dimension $\binom{
    3+4-1}{4} =15$. One computes that $\dim P(\Delta)_{E,3} = 8$.

  Adding all the terms up with the appropriate powers of $x-1$ and
  $y$, Theorem~\ref{thm:main} says that
  \[
  (x-1)^3 + 7 (x-1)^2 + 14(x-1) + 7(x-1)y + y^4 + 3y^3 + 6y^2 + 10y +
  8.
  \]
  is the Tutte polynomial of the Fano matroid.
\end{example}

\section{Hilbert Series}\label{sec:cor}
Recall that the polynomial algebra $\Sym(V^*)$ has a grading
\[
\bigoplus_{d \geq 0} \Sym^d(V^*)
\]
and that a vector subspace $P \subset \Sym(V^*)$ is said to be graded
if it is equal to the direct sum of its homogeneous pieces $P \cap
\Sym^d(V^*)$. The Hilbert series of $P$, $\Hilb(P,t)$, is the
generating function 
\[\sum_{k \geq 0} \dim(P \cap \Sym^k(V^*))\ t^k.\]

Denote the Tutte polynomial of $M(\Delta)$ by $T(\Delta;x,y)$. Since
$P(\Delta)$ is generated by products of linear forms it is a graded
subspace of $\Sym(V^*)$.
\begin{corollary}[Wagner {\cite[Proposition 3.1]{wagner}}]\label{cor:hilb}
  We have,
  \[
  \Hilb(P(\Delta),t) = t^{n-l}T(\Delta;1+t,1/t).
  \]
  The dimension of $P(\Delta)$ is the number of independent sets of
  the matroid $M(\Delta)$.
\end{corollary}
\begin{proof}
  The first claim is found by setting $x=1+t$ and $y=1/t$ in
  Theorem~\ref{thm:main}. The second claim follows from the well-known Tutte
  polynomial evaluation of $T(M;2,1)$ as the number of independent
  sets of the matroid $M$.
\end{proof}
A set $C \subset E$ is called a \textit{circuit} of $M(\Delta)$ if $C$
is dependent but $C-c$ is independent for all $c \in C$. If $I \subset
E$ is independent, define $e \in E = \{1,2,\dots,n\}$ to be
\textit{externally active} in $I$ if $e$ is the smallest element of a
circuit in $I \cup e$. Denote the set of elements externally active in
$I$ by $ex(I)$.
\begin{corollary}[Ardila-Postnikov {\cite[Theorem
    4.2.2]{ardila}}]\label{cor:x=0}
  The Hilbert series of $P(\Delta)_E$, the subspace of $P(\Delta)$
  spanned by those $\alpha_S$ with $\rank(E-S)=\ell$, is
  \[
  \Hilb( P(\Delta)_E ,t) = t^{n-\ell}T(\Delta;1,1/t).
  \]
  The dimension of $\dim P(\Delta)_{X,k}$ is the number of independent
  sets $I$ of $M(\Delta)$ with closure $X$ such that $|I \cup
  ex(I)|=k$.
\end{corollary}

We will need the notion o
\begin{proof}
  Set $x=1$ and $y=1/t$ in Theorem~\ref{thm:main} to obtain the
  Hilbert series. When $X=E$ the second claim follows from the well
  known expression for the Tutte polynomial of a matroid $M$ on an
  ordered set $E$ as
  \[
  T(M;x,y) = \sum_{B \textup{ a base of } M} x^{|in(B)|}y^{|ex(B)|},
  \]
  Here $in(B)$ is the set of elements of $e \in B$ that are the
  smallest elements of a bond of $E-(B-e)$; these are the internally
  active elements of $B$.

  In case $X \neq E$ write $\Delta_X$ for the sequence obtained from
  $\Delta$ by deleting the $\alpha_i$ with $i \notin X$. It follows
  that $P(\Delta)_{X,k}$ is isomorphic to $P(\Delta_X)_{X,k}$, the
  isomorphism being multiplication by $\alpha_{E-X}$. In this way we
  reduce to the case when $X = E$.
\end{proof}
A circuit $C$ with its smallest element deleted is called a broken
circuit. A set $S \subset E$ is said to have \textit{no broken
  circuits}, or be \textit{nbc}, if it does not contain any broken
circuits. This implies that $S$ does not contain any circuit, and
hence, is independent. In terms of external activity, $I$ is nbc if it
is independent and $ex(I) = \emptyset$.
\begin{corollary}[Orlik-Terao {\cite[Theorem
    4.3]{orlik-terao}}]\label{cor:y=0}
  The Hilbert series of 
  \[
  \bigoplus_{X \in L(\Delta)} P(\Delta)_{X,\rank(X)}
  \]
  is the equal to $t^{n-\ell} T(\Delta;1+t,0)$. The dimension of
  $P(\Delta)_{X,\rank(X)}$ is the number of nbc sets of $M(\Delta)$
  with closure $X$.
\end{corollary}
In the case that $\Delta$ is the collection of linear forms defining a
complex hyperplane arrangement $\mathcal{A}$, there is a natural
isomorphism of graded vector spaces between this vector space and the
complexified cohomology of the complement $V-\bigcup_i
\ker(\alpha_i)$. This map takes $\alpha_{E-I}$ to $(\alpha_{i_1}
\wedge \dots \wedge \alpha_{i_k})/\alpha_I$, where $I =
\{i_1<\dots<i_k\}$.
\begin{proof}
  Set $x = 1+1/t$ and $y=0$ in Theorem~\ref{thm:main} to obtain the Hilbert
  series. The second claim follows from Corollary~\ref{cor:x=0}, since
  the dimension of $P(\Delta)_{X,\rank(X)}$ is the number of bases $I$
  of $X$ with $ex(I) = \emptyset$.
\end{proof}

\section{Algebras Generated by Recipricals of Linear Forms}\label{sec:alg}
Let $\Sym(V^*)_{(0)}$ be the $K$-algebra of rational functions on
$V$. Assume that $\alpha_i \neq 0$ for all $i$. In this section we
will investigate the $K$-algebra of non-constant rational functions
without zeros whose poles are contained in the hyperplane arrangement
$\bigcup_{i \in E} \ker \alpha_i$. This was studied eariler by Brion
and Vergne \cite{brion-vergne}, and Terao \cite{terao}.

Let $C(\Delta)$ be the $K$-algebra generated by the rational functions
$\{1/\alpha_i:i \in E\}$. Define the degree of $1/\alpha_i$ to be $1$,
so that $C(\Delta)$ is a graded algebra. If $X \in L(\Delta)$, define
$C(\Delta)_{X}$ to be the space spanned by products $(\alpha_1^{d_1}
\cdots \alpha_n^{d_n})^{-1}$ where the set $\{i: d_i > 0\} \subset E$
has closure $X$ in $M(\Delta)$. Let the degree $k$ piece of
$C(\Delta)_X$ be $C(\Delta)_{X,k}$. Let $m\Delta$ be the sequence
$\Delta$ repeated $m$ times. For example, if $\Delta = (\alpha,\beta)$
then
\[
4 \Delta = (\alpha,\beta,\alpha,\beta,\alpha,\beta,\alpha,\beta).
\]
Since adding parallel elements does not change the span of a set, we
can identify the flats of $m\Delta$ with those of $\Delta$. 
\begin{proposition}\label{prop:alg}
  For each $k \geq 1$ there is an isomorphism of vector spaces
  \[
  C(\Delta)_{X,k} \to P(k\Delta)_{X,k}.
  \]
\end{proposition}
\begin{proof}
  The isomorphism is simply multiplication by $(\alpha_1 \alpha_2
  \cdots \alpha_n)^k$.
\end{proof}
Proposition~\ref{prop:directsum} allows us to conclude the following
result of Terao.
\begin{proposition}[Terao {\cite[Proposition~2.1]{terao}}]
  There is a $K$-vector space direct sum decomposition
  \[
  C(\Delta) = \bigoplus_{\substack{X \in L(\Delta)\\ k \geq 0}}
  C(\Delta)_{X,k}
  \]
\end{proposition}

In \cite{brion-vergne}, Brion and Vergne first studied $C(\Delta)$ and
its subalgebra $C(\Delta)_E$. They viewed the latter as a module for
the algebra $\partial(V)$ of constant coefficient derivations on
$V$. Assuming that $\operatorname{char}(K) =0$, one of their main
results was that $C(\Delta)_E$ is a free $\partial(V)$-module.  In
\cite{terao} Terao, using Brion and Vergne's result, derived a formula
for the Hilbert series of $C(\Delta)$ in terms of the Poincar\'e
polynomial of $\Delta$. The goal of this section is to derive Terao's
formula with no assumption on the characteristic of the field $K$.
\begin{theorem}[Terao {\cite[Theorem 1.2]{terao} if
    $\operatorname{char}(K) =0$}]\label{thm:alg}
  We have,
  \[
  \Hilb(C(\Delta),t) = \left(\frac{y}{1-y}\right)^{\ell}
  T(\Delta;1/y,0).
  \]
\end{theorem}
To prove this we will need the the following easy result, which
appears as Lemma~{6.3.24} in \cite{tutte-poly}. It allows us to
determine the Tutte polynomial of $m\Delta$ in terms of the Tutte
polynomial of $\Delta$.
\begin{lemma}\label{lem:parallel}
  Let $M$ be a matroid and $m$ be a positive integer. If $m M$ is
  the matroid obtained from $M$ by replacing every element of $M$ by
  $m$ parallel elements then
  \[
  T(m M;x,y) = \left(\frac{1-y^m}{1-y}\right)^{\rank(M)}
  T\left(M;\frac{xy-x-y+y^m}{y^m-1},y^m\right).
  \]
\end{lemma}
\begin{proof}[Proof of Theorem~\ref{thm:alg}] For a flat $X \in
  L(\Delta)$, let $\Delta_X$ be the sequence obtained from $\Delta$ by
  deleting those $\alpha_i$ where $i \notin X$. By
  Proposition~\ref{prop:alg} the dimension of $C(\Delta)_{X,k}$ is the
  dimension of $P(k\Delta)_{X,k}$. Multiplication by
  $(\alpha_{E-X})^k$ gives rise to a $K$-vector space isomorphism
  \[
  P(k \Delta_X)_{X,k} \to P(k \Delta)_{X,k}
  \] 
  and so it suffices to compute $\dim P(k\Delta_X)_{X,k}$ to find
  $\dim C(\Delta)_{X,k}$.

  If $f$ is a polynomial in $y$ then $[y^j]f$ is the coefficient of
  $y^j$ in $f$. By Theorem~\ref{thm:main} and Lemma~\ref{lem:parallel} the
  dimension of $P(kX)_{X,k}$ can be written as
  \[ 
  \begin{split}
    [y^{k-\rank(X)}]T(k\Delta_X;1,y)=&
    [y^k]\left(y\frac{1-y^k}{1-y}\right)^{\rank(X)}T(\Delta_X;1,y^k)\\
    =&[y^k]\left(\frac{y}{1-y}\right)^{\rank(X)}T(\Delta_X;1,0)
  \end{split}
  \]
  The second equality follows since the powers of $y$ that appear in
  $T(\Delta_X;1,y^k)$ are multiples of $k$. We know that
  $T(\Delta_X;1,0)$ is the number of nbc bases of $X$, that is, the
  number of nbc sets of $M(\Delta)$ with closure $X$, and it follows
  that
  \[
  \Hilb(C(\Delta)_X,t) = |\{\textup{nbc bases of }X\}|
  \left(\frac{t}{1-t}\right)^{\rank(X)}.
  \]
  Summing over all flats $X \in L(\Delta)$,
  \[
  \begin{split}
    \Hilb(C(\Delta),y) =& \sum_{X \in L(\Delta)} |\{\textup{nbc bases
      of } X\}| \left(\frac{y}{1-y}\right)^{\rank(X)}\\
    = & \left(\frac{y}{1-y}\right)^\ell T(1/y,0).
  \end{split}
  \]
  The second equality can be verified using, e.g., the Tutte
  polynomial interpretation in Corollary~\ref{cor:y=0}.
\end{proof}
\section{Spanning Homogeneous Pieces of Symmetric Powers}\label{sec:ff}
In this section we investigate the following problem: Given a sequence
of linear forms $\Delta = (\alpha_1,\dots,\alpha_n)$ spanning the dual
of a $K$-vector space $V$, determine those $d$ such that every $d$-form
on $V$ can be written in the form
\[
\sum_{S \in \binom{E}{d}} c_S \alpha_S.
\]
Phrased differently, determine those $d$ such that $\Sym^d(V^*)
\subset P(\Delta)$. The answer to this question is phrased in terms of
the cociruits. Recall that collection of \textit{cocircuits} of
$\Delta$ is
\[
\{E-H: H \in L(\Delta), \rank(H) = \ell-1\}.
\]
It is possible to reconstruct the matroid $M(\Delta)$ from its
collection of cocircuits.
\begin{theorem}\label{thm:spanning}
  There is a containment $\Sym^d(V^*) \subset P(\Delta)$ if and only
  if $d$ is less than or equal to the size of the smallest cocircuit
  of $\Delta$.
\end{theorem}

\begin{example}
  For $q$ a prime power, consider all of the $[\ell]_q :=
  (q^\ell-1)/(q-1)$ hyperplanes in $V=\mathbb{F}_q^\ell$. Let $\Delta$
  be a sequence of linear forms, one defining each of these
  hyperplanes. It is a fact that every cocircuit of $\Delta$ has size
  $[\ell]_q - [\ell-1]_q = q^{\ell-1}$. Theorem~\ref{thm:spanning}
  states that $P(\Delta)$ contains $\Sym^d(V^*)$ if and only if $d
  \leq q^{\ell-1}$.

  In the case that $q=2$ and $\ell=3$ then $\Delta$ is the collection
  from Example~\ref{ex:fano}. Using the Tutte polynomial calculation
  there,
  \[
  \Hilb(P(\Delta),t) = 1+3t+6t^2+10t^3+15t^4+14t^5+7t^6+t^7.
  \]
  For $d \leq 2^{3-1} = 4$ the coefficient of $t^d$ can be written as
  $\binom{3+d-1}{d}$. For these $d$, $P(\Delta) \cap \Sym^d(V^*)$ and
  $\Sym^d(V^*)$ have equal dimensions, so they are equal.
\end{example}
\begin{proof}[Proof of Theorem~\ref{thm:spanning}]
  To start, we investigate the Hilbert series of $P(\Delta)$ and
  determine the first coefficient not of the form
  $\binom{\ell+d-1}{d}$. By Corollary~\ref{cor:hilb},
  \begin{align*}
    \Hilb(P(\Delta),t)=
    t^{n-\ell} T\left(\Delta;1+t,t^{-1}\right)
  \end{align*}
  Rewrite this using the formula,
  \[
  T(M;x,y) = \frac{1}{(y-1)^{\rank(M)} }\sum_{X \in L(M)} y^{|X|}
  \chi(M/X;(x-1)(y-1),0).
  \]
  Here $\chi(M;\lambda)$ is the \textit{characteristic polynomial} of
  $M$, obtained from $T(M;x,y)$ by the rule $\chi(M;\lambda) =
  (-1)^{\rank(M)}T(1-\lambda,0)$. This formula can be found in the
  discussion of the \textit{coboundary polynomial} in
  \cite{tutte-poly}. Doing the stated substitutions allows us to write
  \[
  \Hilb(P(\Delta),t) = \frac{1}{(1-t)^\ell} \sum_{X \in L(\Delta)}
  t^{n -|X|} (-1)^{\ell-\rank(X)}T(M(\Delta)/X;t,0)
  \]
  Denote the polynomial on the right by $h(t)$. If $h(t)$ is of the
  form $1-at^{N+1} + \cdots$, where $a$ is nonzero and the ellipsis
  denotes higher degree terms, then
  \[
  \Hilb(P(\Delta),t) = \sum_{k=0}^N \binom{\ell+k-1}{k} t^k +
  \left(\binom{\ell+N}{N+1} - a\right)t^{N+1} + \cdots
  \]
  This is to say, $P(\Delta)$ contains $\Sym^d(V^*)$ for $d \leq N$
  and $P(\Delta)$ does not contain all of $\Sym^{N+1}(V^*)$.  To find
  the smallest power of $t$ appearing in $h(t)$ consider the smallest
  power of $t$ appearing in each of its summands. If $X \neq E$ then
  $T(M(\Delta)/X;t,0)$ has no constant term and hence the smallest
  power of $t$ in each summand is at least $n-|X|+1$. The summand
  corresponding to $X=E$ is the constant polynomial $1$. Let $X_0$
  denote a flat of largest size which is not $E$. For any such flat
  $M(\Delta)/X_0$ is a rank one matroid with no loops. This implies
  that $T(M(\Delta)/{X_0};t,0) = t$ and hence
  \[
  h(t) = 1 - a t^{n-|X_0|+1} + \cdots
  \]
  where $a$ is the number of flats of $\Delta$ with size $|X_0|$. To
  summarize: \textit{If $d \leq |E-X_0|$ then $P(\Delta)$ contains
    $\Sym^d(V^*)$. Further, $P(\Delta)$ does not contain all of
    $\Sym^{|E-X_0|+1}(V^*)$.}

  Suppose there was some $d > |E-X_0|+1$ such that $P(\Delta)$
  contained $\Sym^d(V^*)$. If $\alpha_1,\dots,\alpha_\ell$ are a basis
  for $V^*$ then for all sequences $(\sigma(1),\dots,\sigma(\ell))$
  such that $\sigma(1)+\dots+\sigma(\ell) = d$,
  \[
  \alpha_1^{\sigma(1)} \cdots \alpha_\ell^{\sigma(\ell)} \in P(\Delta)
  \]
  We will prove in Lemma~\ref{lem:ex} that this implies
  \[
  \alpha_1^{\tau(1)} \cdots \alpha_\ell^{\tau(\ell)} \in P(\Delta)
  \]
  for any $\tau$ such that $\tau(i) \leq \sigma(i)$. It follows
  that\footnote{If $\operatorname{char}(K)=0$ it is easy to see this
    without the lemma. Indeed, $P(\Delta)$ is a module for
    $\partial(V^*)$, the $K$-algebra of constant coefficient
    differential operators on $V^*$. For $1 \leq i \leq \ell$, there
    is a unique differential operator on $V^*$ taking $\alpha_i$ to
    $1$ and $\{\alpha_1,\dots,\alpha_\ell\}-K\cdot \alpha_i$ to zero.
    Thus, if $\alpha_1^{\sigma(1)} \cdots \alpha_\ell^{\sigma(\ell)}$
    is in $P(\Delta)$ then by applying these differential operators we
    get that an integer multiple of $\alpha_1^{\tau(1)} \cdots
    \alpha_\ell^{\tau(\ell)}$ is in $P(m\Delta)$ whenever $\tau(i)
    \leq \sigma(i)$.} $P(\Delta)$ contains $\Sym^{d'}(V^*)$ for any
  $d' \leq d$. This cannot be, since we known that $P(\Delta)$ does
  not contain $\Sym^{d'}(V^*)$ when $d'= |E-X_0|+1<d$. It follows that
  $P(\Delta)$ contains $\Sym^d(V^*)$ if and only if $d \leq |E-X_0|$.
  By definition, $E-X_0$ is a cocircuit of smallest size, so the
  theorem follows.
\end{proof}


\section{Deletion-Contraction Proof of Theorem~\ref{thm:main}}\label{sec:proof}
Define $H(\Delta;x,y)$ by the rule
\[
H(\Delta;x,y) = \sum_{\substack{ X \in L(\Delta)\\ k \geq 0}}
x^{\rank(M(\Delta))-\rank(X)} y^{k - \rank(X)} \dim P(\Delta)_{X,k}.
\]
Here $\rank(M(\Delta))$ is the rank of the matroid of $M(\Delta)$,
which is $\ell$. Theorem~\ref{thm:main} can be stated as 
\[
H(\Delta;x,y) = T(\Delta;1+x,y).
\]
To prove this we will apply the following fundamental theorem on the
Tutte polynomial (see \cite[Theorem~6.2.2]{tutte-poly}).
\begin{theorem}\label{thm:tutte}
  Let $\mathcal{M}$ be the class of isomorphism classes of
  matroids. There is a unique function, called the Tutte polynomial,
  $T:\mathcal{M} \to \mathbb{Z}[x,y]$ which satisfies
  \begin{enumerate}
  \item[(a)] The Tutte polynomial of the one element isthmus
    is $x$ and that of the one element loop is $y$.
  \item[(b)] The Tutte polynomial is multiplicative in the sense that
    $T(M \oplus N) = T(M)T(N)$, where $M \oplus N$ is the direct sum
    of matroids.
  \item[(c)] If $M$ is a matroid on $E$ and $e \in E$ is
    neither a loop nor an isthmus then
    \[
    T(M) = T(M-e) + T(M/e).
    \]
    Here $M-e$ is deletion of $e$ from $M$ and $M/e$ is contraction of
    $M$ by $e$.
  \end{enumerate}
\end{theorem}
To prove that $H(\Delta;x,y)$ equals $T(\Delta;1+x,y)$ we need to
check that it satisfies properties (a), (b) and (c) in
Theorem~\ref{thm:tutte}.
\begin{proof}[Verification of property (a)]
  It must be checked that
  \[
  \begin{split}
    H(\{1\};x,y) =& T(\{1\};1+x,y)= 1+x,\\
    H(\{0\};x,y) =& T(\{0\};1+x,y)= y.
  \end{split}
  \]
  This easy task is left to the reader.
\end{proof}
\begin{proof}[Verification of property (b)]
  Suppose that 
  \[
  \Delta = (\alpha_1,\dots,\alpha_n), \quad \Delta' =
  (\beta_1,\dots,\beta_m)
  \] are sequences of linear forms on two vector spaces $V$ and $W$
  over $K$. Let $\Delta \oplus \Delta'$ denote the concatenation of
  the sequences, $(\alpha_1,\dots,\alpha_n,\beta_1,\dots,\beta_m)$,
  viewed as vectors in $V^* \oplus W^*$. This agrees with the direct
  sum of matroids since the matroid of $\Delta \oplus \Delta'$ is the
  direct sum of the matroids of $\Delta$ and $\Delta'$. Recall that
  there is a natural isomorphism of graded $K$-algebras,
  \begin{align}\label{natisom}
    \Sym(V^* \oplus W^*) \approx \Sym(V^*) \otimes \Sym(W^*).
  \end{align}
  The flats of $\Delta \oplus \Delta'$ are in bijection with
  $L(\Delta) \times L(\Delta')$. We claim that if $X \in L(\Delta)$
  and $Y \in L(\Delta')$ then, as graded vector spaces,
  \[
  P(\Delta \oplus \Delta')_{(X, Y)} \approx P(\Delta)_X \otimes
  P(\Delta')_Y
  \]
  Indeed, the isomorphism \eqref{natisom} maps $\alpha_S \otimes
  \beta_T$ to $\alpha_S\beta_T$ which is in $P(\Delta \oplus
  \Delta')_{(X,Y)}$. Since every monomial defining $P(\Delta \oplus
  \Delta')_{(X,Y)}$ is of this form, we have the needed
  isomorphism. Lastly, since the rank of the flat corresponding to
  $(X,Y)$ is sum of the ranks of $X$ and $Y$ and the isomorphism
  \eqref{natisom} is of \textit{graded algebras},
  \[
  H(\Delta \oplus \Delta'; x,y) = H(\Delta;x,y)H(\Delta';x,y).\qedhere
  \]
\end{proof}

To set up our verification of property (c), note that for any linear
form $\alpha \in V^*$ there is a complex of graded vector spaces
\begin{align}\label{eq:exseq0}
0 \to \Sym(V^*) \stackrel{\cdot \alpha}{\to} \Sym(V^*)
\to \Sym(\ker(\alpha)) \to 0.
\end{align}
The second map is induced by restricting linear forms on $V$ to
$\ker(\alpha)$. When $\alpha \neq 0$ this complex is exact by
definition. Let $\Delta - \alpha_i$ be obtained from $\Delta$ by
deleting the element in the $i$-th position and $\Delta/\alpha_i$ be
the restriction of the forms in $\Delta-\alpha_i$ to
$\ker(\alpha^*)$. These definitions are made in such a way that they
coincide with deletion and contraction of matroids:
\[
M(\Delta - \alpha_i) = M(\Delta)-i, \quad M(\Delta/\alpha_i) = M/i
\]
where $M(\Delta), M(\Delta-\alpha_i)$ and $M(\Delta/\alpha_i)$
denote, respectively, the matroids of $\Delta$, $\Delta-\alpha_i$
and $\Delta/\alpha_i$. 

For any $\alpha_i$, the complex \eqref{eq:exseq0} is seen to
restrict to the complex
\begin{align}\label{exseq}
  0 \to P(\Delta-\alpha_i) \stackrel{\cdot \alpha_i}{\to} P(\Delta)
  \to P(\Delta/\alpha_i) \to 0.
\end{align}
\begin{lemma}\label{lem:ex}
  If $i \in E$ is neither a loop nor an isthmus of $M(\Delta)$ then
  this complex is exact.
\end{lemma}
The assumption that $i$ is not an isthmus is strictly for technical
reasons: We need to have that $\Delta-\alpha_i$ spans $V^*$, which
will not happen if $i$ is an isthmus.
\begin{proof}
  The existence of the complex guarantees that $$\dim P(\Delta) \geq
  \dim P(\Delta-\alpha_i) + \dim P(\Delta/\alpha_i).$$ We may assume,
  by induction, that if $\Delta'$ has fewer elements than $\Delta$
  then $\dim P(\Delta')$ is the number of independent sets of the
  matroid $M(\Delta')$. It follows, from our induction hypothesis,
  that $\dim P(\Delta)$ is at least the sum of the number of
  independent subsets of $M(\Delta-\alpha_i)$ and
  $M(\Delta/\alpha_i)$. This is the number of independent sets of
  $M(\Delta)$, so we need only show that $P(\Delta)$ is spanned by
  this many elements to conclude exactness.

  We claim that the elements $\{p_{E-(I \cup ex(I))}: I \textup{ is
    independent in } M(\Delta)\}$ span $P(\Delta)$, where $ex(I)$ is
  the set of elements which are externally active in $I$ (\textit{cf.}
  Section~\ref{sec:cor}). The proof of the claim is by what Las
  Vergnas and Forge call \textit{lexicographic compression} and
  variants of this proof abound (see e.g., Ardila's thesis
  \cite{ardila} a nearly identical claim. This claim also appears in
  \cite{ardila-postnikov}.). Let $S$ be the lexicographically least
  subset of $E$ such that $\alpha_{E-S}$ is not in the span of these
  elements. We can uniquely write $S = I \cup J$ where $J \subset
  ex(I)$. To see this we take $I \subset S$ to be the
  lexicographically largest independent set spanning $S$, and $J$ to
  be the complement of $I$ in $S$. One immediately sees that $J
  \subset ex(I)$. We will be done if we can show that $J = ex(I)$.
  Suppose, on the contrary, that $f \in ex(I)-J$ and write $\alpha_f=
  \sum_{e \in I} c_e \alpha_e$. We have $\alpha_{E-S} = \alpha_{E-(S
    \cup f)} \alpha_f$ and hence
  \[
  \alpha_{E-S} = \sum_{e \in I} c_e \alpha_{E-(I \cup J \cup f)}
  \alpha_e = \sum_{e \in I} c_e \alpha_{E-((I-e) \cup (J \cup f))}
  \]
  Since $\alpha_{E-S}$ is not in the span of $\{ p_{E -(I \cup ex(I))}
  : I \in \mathcal{I}(M(\Delta))\}$ we know there is some element on
  the right which is not in the span of these elements. For each $e
  \in I$ we see that $(I-e) \cup (J \cup f)$ is lexicographically
  smaller than $I \cup J$, which is a contradiction.
\end{proof}
\begin{remark}
  The lemma proves Corollary~\ref{cor:x=0} by exhibiting a basis for
  $P(\Delta)_{X,k}$ with the cardinality stated there.
\end{remark}
There are two ways to proceed from here. The first is to use the
combinatorial basis just constructed to verify that the
deletion-contraction recurrence holds. The second is show that the
exact sequence above can be refined to consider a flat and a
degree. In the linear-algebraic spirit of this paper, we take the
latter route.
\begin{proof}[Verification of property (c)]
  Assume that $i \in E$ is neither a loop nor an isthmus of
  $M(\Delta)$, that is, $r(i) = 1$ and $\rank(E-i) = \rank(E)$.
  Recall (see \cite[Chapter~7]{matroid-th}) that the flats of the
  deletion $\Delta-\alpha_i$ are in bijection with the flats of $X$ of
  $\Delta$ such that $\rank(X-i)=\rank(X)$. Also, the flats of
  $\Delta/\alpha_i$ are bijection with the flats of $\Delta$
  containing $i$. We wish to refine the exact sequence \eqref{exseq}
  to consider a flat $X$ of $\Delta$ and a degree. To do this,
  consider three cases depending on whether or not $i \in X$, and if
  $i \in X$, then whether $i$ is an isthmus of $X$.
  
  First suppose that $i \notin X$. It follows that $X$ is a flat of
  $\Delta-\alpha_i$. The complex \eqref{exseq} restricts to the exact
  complex
  \begin{align}\label{eq:delcontr1}
    0 \to P(\Delta-\alpha_i)_{X,k} \stackrel{\cdot \alpha_i}{\to}
    P(\Delta)_{X,k} {\to} 0.
  \end{align}
  Every product $\alpha_S \in P(\Delta)_{X,k}$ is of the form
  $\alpha_S$ where $i \in S$, hence $\alpha_{S-i}$ lies in
  $P(\Delta-\alpha_i)_{X,k}$ and the map is surjective. Since the map
  is the restriction of an injection it is an isomorphism.
  
  Suppose that $i \in X$ and $i$ is not an isthmus of $X$. In this
  case $X-i$ is a flat of both $\Delta-\alpha_i$ and $\Delta/\alpha_i$
  and we claim that \eqref{exseq} restricts to the exact complex
  \begin{align}\label{eq:delcontr2}
    0 \to P(\Delta-\alpha_i)_{X-i,k} \stackrel{\cdot \alpha_i}{\to}
    P(\Delta)_{X,k} \to P(\Delta/\alpha_i)_{X-i,k-1} \to 0.
  \end{align}
  To see this we pick some $\alpha_S \in P(\Delta-\alpha_i)_{X-i,k}$
  and see that $\alpha_{S \cup i}$ is in $P(\Delta)_{X,k}$. This is
  because $(E-i)-S = E-(S\cup i)$ has closure $X$ in $M(\Delta)$. If
  the closure were the smaller set $X-i$, then $i$ would be have been
  an isthmus of the flat $X$. The map on the left in
  \eqref{eq:delcontr2} is the restriction of an injection, hence we
  have exactness on the left.
  
  If $\alpha_S \in P(\Delta)_{X,k}$ and $i \notin S$ then the closure
  in $M(\Delta/\alpha_i)$ of $(E-i)-S$ is $X-i$. The degree of
  $\alpha_S$ is unchanged under $\alpha_i\mapsto 0$ hence
  $P(\Delta)_{X,k}$ has image in $P(\Delta/\alpha_i)_{X-i,k-1}$. That
  every monomial spanning $P(\Delta/\alpha_i)_{X-i,k-1}$ is in the
  image of $P(\Delta)_{X,k}$ follows from the definition of
  $\Delta/\alpha_i$ as the restriction of the elements of
  $\Delta-\alpha_i$ to $\ker(\alpha_i)$. The exactness on the right of
  \eqref{eq:delcontr2} follows.

  We now prove exactness in the middle of \eqref{eq:delcontr2}. If an
  element of $P(\Delta)_{X,k}$ restricts to zero on $\ker(\alpha_i)$
  then, by Lemma~\ref{lem:ex}, it can be written as $\alpha_i$ times
  some linear combination $\sum c_S \alpha_S$ where $\alpha_S \in
  P(\Delta-\alpha_i)$. For these $S$ we have $\alpha_i \alpha_S \in
  P(\Delta)_{X,k}$ and, since $i$ was not an isthmus of $X$, we have
  $\alpha_S$ is in $P(\Delta-\alpha_i)_{X-i,k}$.
 
  In the case that $i \in X$ and $i$ is an isthmus of $X$ it follows
  that \eqref{exseq} restricts to the exact complex
  \begin{align}\label{eq:delcontr3}
    0 \to P(\Delta)_{X,k} \to P(\Delta/\alpha_i)_{X-i,k-1} \to 0.
  \end{align}
  The surjectivity here is clear. If an element is in the kernel of
  this map we may, as before, write it as a linear combination of
  terms $\alpha_i\alpha_S \in P(\Delta)_{X,k}$. It follows that $E-(S
  \cup i)$ has closure $X$ in $M(\Delta)$, which cannot be since $i
  \notin E-(S\cup i)$ and $i$ is in every base of $X$. We conclude
  that the kernel is zero and \eqref{eq:delcontr3} is exact.

  Finally, we can verify the deletion-contraction recurrence. To do
  so, break up the defining sum for $H(\Delta;x,y)$ according to the
  three cases we just considered. Let $L_1 \subset L(\Delta)$ be of
  the set of flats of $\Delta$ not containing $i$, $L_2 \subset
  L(\Delta)$ be the set of flats containing $i$ as an isthmus, and let
  $L_3 \subset L(\Delta)$ be the remaining flats. If $\corank(X)$
  denotes $\rank(M(\Delta))-\rank(X)$, the corank of $X$, we see that
  the exact complexes \eqref{eq:delcontr1}, \eqref{eq:delcontr2} and
  \eqref{eq:delcontr3} imply that we can write
  \begin{align*}
    H(\Delta;x,y) =& \sum_{X \in L_1 \cup L_3,k}x^{\corank(X)} y^{k -
      \rank(X)}\dim
    P(\Delta-\alpha_i)_{X-i,k}\\
    & \quad + \sum_{X \in L_2 \cup L_3,k} x^{\corank(X)} y^{k -
      \rank(X)} \dim P(\Delta/\alpha_i)_{X-i,k-1}.
  \end{align*}
  The first sum is \eqref{eq:delcontr1} and the left of
  \eqref{eq:delcontr3}, while the second sum is \eqref{eq:delcontr2}
  and the right of \eqref{eq:delcontr3}. We claim that the first sum
  here is $H(\Delta-\alpha_i)$ and the second is
  $H(\Delta/\alpha_i)$. Since the flats of $\Delta-\alpha_i$ are in
  bijection with $\{X-i: X \in L_1 \cup L_3\}$ and the flats of
  $\Delta/\alpha_i$ are in bijection with $\{X-i: X \in L_2 \cup
  L_3\}$, we only need to check that the exponents of $x$ and $y$ in
  each sum are correct. If $X \in L_1 \cup L_3$, then the rank of
  $X-i$ in the matroid $M(\Delta-\alpha_i)$ is equal to the rank of
  $X$ in $M(\Delta)$.  If $X \in L_2 \cup L_3$ then the rank of $X-i$
  in $M(\Delta/\alpha_i)$ is one less than its rank in $M(\Delta)$. It
  follows that the exponents of $x$ and $y$ are correct and so
  \begin{align*}
    H(\Delta;x,y) = H(\Delta-\alpha_i;x,y) + H(\Delta/\alpha_i;x,y),
  \end{align*}
  which is what we wanted to show.
\end{proof}

\section{A Second Proof of Theorem~\ref{thm:main}}\label{sec:proof2}
It has been observed by an anonymous referee that
Theorem~\ref{thm:main} may be proved directly from Ardila and
Postnikov's result in Corollary~\ref{cor:x=0}. This is the case, and
in this section we give this proof. Note that the proof of
Corollary~\ref{cor:x=0} in \cite{ardila} is characteristic
independent, but the proof in \cite{ardila-postnikov} is not.

Recall that Corollary~\ref{cor:x=0} states, among other things, that
\[
\dim P(\Delta)_{X,k}
\]
is the number of independent sets $I$ of $M(\Delta)$ with closure $X$
such that $|I \cup ex(I)| = k$. 

Also, the Tutte polynomial of an arbitrary matroid $M$ can be written
as
\[
T(M;x,y) = \sum_B x^{|in(B)|} y^{|ex(B)|}
\]
the sum over bases $B$ of $M$. Here $in(B)$ is the set elements that
are internally active in $B$ and $ex(B)$ is the set of elements that
are externally active in $B$. Recall that $e \in B$ is internally
active in $B$ if $e$ is the smallest element of a bond of $D \subset
E-(B-e)$, i.e., $e$ is the smallest element $f$ of a (necessarily
unique) set $D \subset E-(B-e)$ that is minimal with the property that
\[
r(E-D)< r(E).
\]
It follows from that we may write
\[
T(M;1+x,y) = \sum_B \sum_{I \subset in(B)} x^{|I|} y^{ex(B)}.
\]
We now need a combinatorial result of Crapo.
\begin{lemma}[Crapo {\cite[Lemmas~6,8,9]{crapo}}]
  Let $M$ be an arbitrary matroid on an ordered set $E$. Every subset
  $S \subset E$ can be written uniquely as $S = B-I \cup J$ where $I
  \subset in(I)$ and $J \subset ex(B)$.

  Further, if $I \subset in(B)$ then $ex(B-I) = ex(B)$.
\end{lemma}
It follows that given an independent set $I$, we may write this set
uniquely as $I = B-I_0$, where $I_0 \subset in(B)$. From this, we may
write the Tutte polynomial as
\[
\begin{split}
  T(M;1+x,y) &= \sum_B \sum_{I \subset in(B)} x^{|B-I|
  }y^{|ex(B-I)|}  \\
  &= \sum_I x^{r(M)-r(I)} y^{|ex(I)|}
\end{split}
\]
the second sum over all independent sets $I$ of $M$.
From this it follows at once that
\[
T(M;1+x,y) = \sum_{\substack{X \in L(M)\\k \geq 0}} x^{r(M)-r(X)}
y^{k-r(X)}\cdot \#\left\{I \subset E:
\begin{matrix}
  I \textup{ independent,}\\
  cl(I) = X,\\
  |I \cup ex(I)| = k.
\end{matrix}
 \right\}
\]
Theorem~\ref{thm:main} now follows from this expression for
$T(M;1+x,y)$ and Corollary~\ref{cor:x=0}.
\section{Acknowledgements}
The author thanks Victor Reiner for many helpful conversations and
suggesting the problem considered in Section~\ref{sec:ff}. Thanks to
Hiroaki Terao for prompting me to go back to his paper \cite{terao},
and several anonymous referees for their many helpful suggestions. The
author was partially supported through NSF grant DMS 0601010.

\bibliography{tutte}\nocite{*} \bibliographystyle{plain}
\end{document}